\documentclass[twoside, 10pt]{article}
 \RequirePackage{amsgen}
\RequirePackage{amsmath} \RequirePackage{amstext}
\RequirePackage{amsbsy} \RequirePackage{amsopn}
\RequirePackage{amsthm} \RequirePackage{amssymb}
\RequirePackage{epsfig}

\usepackage[mathscr]{eucal}

\theoremstyle{plain}
\newtheorem{theorem}{Theorem}[section]
\newtheorem{proposition}[theorem]{Proposition}
\newtheorem{lemma}[theorem]{Lemma}

\newtheorem{corollary}[theorem]{Corollary}

\let\a\alpha      
  \let\n\nu \let\w\omega
  
  \let\l\lambda   
      
\let\L\Lambda

\def\sD{\sf D}

\def\F{{\mathbb F}}
\def\Fr{{\F_{r}}}
\def\Fq{{\F_{q}}}

\def\Fp{{\F_{p}}}
\def\Z{{\mathbb Z}}
\def\Zp{{\mathbb{Z}_{p}}}
\def\Q{{\mathbb Q}}
\def\Qp{{\Q_p}}

\begin{document}

\pagestyle{myheadings}

\markboth{Sangtae Jeong}{On a question of Goss}

\title{{\large \bf  On a question of Goss}}

\author{ Sangtae Jeong}

\footnotetext{{\em Keywords}: Locally analytic endomorphisms; 1-units \\
This work was  supported by Inha University Research Grant\\
{\em Mathematics Subject Classification}  2000: 11K41.}

\maketitle

\bibliographystyle{alpha}

\begin{abstract}
In this note we answer the question raised by D. Goss in
[Applications of non-Archimedean integration to the $L$-series of
$\tau$-sheaves, {\em J. Number Theory,} 110 (2005), no. 1,
83--113] by proving that the group of locally analytic
endomorphisms on the 1-units of a locally compact field of
characteristic $p>0$ is isomorphic to the $p$-adic integers.

\end{abstract}

\vskip 10pt
 After his proof of the analytic continuation of
characteristic $p$ valued  $L$-function, D. Goss \cite{G2}
remarked that the proof depends crucially on locally analyticity
of exponential functions $u^y$ with $p$-adic integers $y$ on the
1-units of a locally compact field of characteristic $p>0.$ He
also pointed out that his proof would automatically work for any
other type of locally analytic endomorphism on the group of
1-units, if such an endomorphism exists. Then the following
question was raised by Goss, who thought it to reasonable to
expect a negative answer:

\noindent Question: Let $\Fr$ be a finite field of $r$ elements
and $r$ be a power of a prime integer $p.$ Let $\Fr((1/T))$ be the
completion of a rational function field $\Fr(T)$ at the infinite
prime and  $U_1$ be the group of 1-units of $\Fr((1/T)).$ Does
there exist a locally analytic endomorphism $f: U_1 \rightarrow
U_1$ which is not of the form $u \mapsto u^y$ for some $p$-adic
integer $y$ ?

The purpose of this note is to give a negative answer to the
Question above in  very great generality. Indeed, we use the Hasse
derivatives to give a proof that the group of locally analytic
endomorphisms on the  group of 1-units of a locally compact field
of characteristic $p>0$ is isomorphic to the $p$-adic integers. As
a consequence, we show that the group of locally analytic
automorphisms on the group of 1-units of a locally compact field
of characteristic $p>0$ is isomorphic to the $p$-adic units.

\section {Notations and statements of results}
Let $K$ be a non-Archimedian, locally compact local field of
characteristic $p>0,$ with maximal compact subring of $R_K$ and
associated maximal ideal $M_K.$ Let $\F$ be the field of constants
of $K$ and $\pi$ be a prime element of $K.$ Then $K$ can be
identified with the field of the Laurent series in $\pi$ over the
finite field $\F:$
$$K = \F((\pi)).$$
We fix  an absolute value $|\cdot|$ associated to the additive
valuation $v$ on  $K$ so that $|x| = q^{-v(x)}$ where $q$ denotes
the order of $\F.$ Let $\rho \in R_K$ with  $0<t=|\rho|$ and  $\a$
be another element of $R_K.$ The closed ball $B_{{\a}, t}$ around
$\a$ of radius $t$  is  defined  by
$$B_{\a,t}= \{ u \in R_K | \quad |u-{\a}| \leq t \}.$$
Let $U:=U_1$ denote the group of 1-units of $K:$
$$U= 1+ M_K=1+ \pi \F[[\pi]].$$
Then $U \subset R_K$ is the closed ball around 1 of radius $1/q,$
i.e., $U=B_{1,1/q}.$

We say that a continuous function $f:U \rightarrow K $ is analytic
on $U$ if
 and only if $f$ may be expressed as the Tayler series of the form :
\begin{eqnarray}\label{Taylerexp0}
 f(u) = \sum_{n=0}^{\infty}b_n(\frac{u-1}{\pi})^n,
\end{eqnarray}
where $\{b_n\}_{n\geq0} \subset K$ and $b_{n} \rightarrow 0$ as $n
\rightarrow \infty.$

A continuous function $f:U \rightarrow K $ is said to be locally
analytic on $U$ if for each $\a \in U$ there exists $t_{\a} >0$
such that $f$ is analytic on $B_{{\a}, t_{\a}}$ in $U.$

Hence, an analytic function $f$ on $B_{{\a}, t_{\a}}$ can be
expressed as the Tayler series of the form :
\begin{eqnarray}\label{Taylerexp0}
 f(u) = \sum_{n=0}^{\infty}c_n(\frac{u-\a}{\rho})^n,
\end{eqnarray}
where $\{c_n\}_{n\geq0}$ is a null sequence  in $ K$ and $t_\a
=|\rho|.$

Throughout, for $u \in U$ we set $u=1+x$ with $x \in M_K=\pi
\F[[\pi]]$ and  view a locally analytic function $f$ on $U$ as a
function of $x$ but not $u,$ unless otherwise specified. By a
change of variable we see from Equation (\ref{Taylerexp0}) that a
locally analytic function $f$ on $U$  may be rewritten  as
\begin{eqnarray}\label{x-exp}
f(u) =f(1+x)= \sum_{n=0}^{\infty}a_nx^n, \quad {\rm for \quad
some} \quad \{ a_n\}_{n\geq0} \quad{\rm in}\quad  K .
\end{eqnarray}

We remark that if the sequence of coefficients $\{a_n\}_{n\geq0}$
in $K$ is bounded, the series in Equation (\ref{x-exp}) converges
on $M_K$ since $a_nx^n \rightarrow0$ as $n \rightarrow \infty.$ In
what follows we denote by $\Lambda$ the set of locally analytic
endomorphisms $f$ from $U$ to itself:
$$ \Lambda = \{ f: U \rightarrow U | \quad f \quad {\rm is\quad locally \quad analytic\quad
on\quad} U \quad{\rm and} \quad f(uv)=f(u)f(v) {\rm \quad for\quad
all \quad} u , v \in U \}.$$ Then it is obvious that $\L$ has two
binary operations, usual multiplication and composition. We also
see easily that $\L$ is a group under multiplication but not under
composition. On the other hand, the $p$-adic integers $\Zp$ act
naturally on $U$(see \cite{G1,H2}), so that the set $\{u^y : y \in
\Zp\}$ is a subgroup of $ (\L,\cdot)$ since $u^y$ is analytic on
$U,$ hence locally analytic on $U.$  We here refer the reader to
\cite{Am,G2} for details on locally analytic functions.

The following theorem answers the Goss's question:
\begin{theorem}\label{main}
Let $\L$ be the group of locally analytic endomorphisms on $U.$
Then $(\L,\cdot)$ is isomorphic to $(\Zp, +).$
\end{theorem}

We might be wondering what happens if endomorphisms in $\L$ are
replaced with automorphisms. The next result answers this
question:

\begin{corollary}\label{auto}
Let $\L_{0}$ be the group of locally analytic automorphisms on
$U.$ Then $(\L_{0},\circ)$ is isomorphic to $(\Zp^*, \cdot).$
\end{corollary}

We have an interesting result on ultimately periodic sequences
which arise from $p$-adic binomial coefficients. By an ultimately
periodic sequence $\{s_n \}_{n \geq 0}$ we mean here that it is
periodic from some index $n$ on i.e., there exist integers $r>0$
and $\w\geq0$ such that $s_{n+r}=s_n$ for all $n\geq \w.$


\begin{corollary}
The following are equivalent: For $y \in \Zp,$

\noindent (1) $(1+x)^y =\sum_{n=0}^{\infty}\binom{y}{n}x^n \in
\Fp[[x]]$ is rational in $x \in M_K.$

\noindent (2) $\{\binom{y}{n} \}_{n \geq 0}$  is  an ultimately
periodic sequence in $\Fp.$

\noindent (3) $y$ is in $\Z.$

\end{corollary}

\section{Proofs}
We begin by recalling that $(\L,\cdot)$ is the group of locally
analytic endomorphisms on $U.$ We first see that any locally
analytic function $f \in \L$ has expansion coefficients which is
bounded. Indeed, the coefficients belong to the prime subfield
$\Fp$ of $K$ in the following lemma.

\begin{lemma}\label{fpele}
If $f(1+x)=\sum_{n=0}^{\infty}a_nx^n \in 1+xK[[x]]$ is a locally
analytic endomorphism on $U$ in $\L$ Then $a_n \in \Fp$ for all
$n=0,1, \cdots.$
\end{lemma}
\begin{proof}
Since a $p$th power mapping on $U$ is linear, and a locally
analytic function $f \in \L$ is an endomorphism on $U$ we obtain
$f(1+x^p)= f((1+x)^p)=f(x+1)^p.$ By equating expansion
coefficients of two end series we have $ a_n^p =a_n$ for all
$n=0,1, \cdots.$ Hence the result follows.

\end{proof}

\begin{lemma}\label{ancoeff}
Let $f(1+x)=\sum_{n=0}^{\infty}a_nx^n \in 1+x\Fp[[x]]$ be a
locally analytic endomorphism in $\L.$ Then $a_n$ is zero if and
only if ${\sD}^{(n)}f(1+x)$ is identically zero, where
${\sD}^{(n)}f$ denotes the Hasse derivative of order $n$ of $f.$
In particular, $a_1$ is zero if and only if $f(1+x) =g(1+x)^p$ for
some $g(1+x)$ in $\L.$
\end{lemma}
\begin{proof}
Since a locally analytic function $f \in \L$ is an endomorphism on
$U$ we can write, for two variables $x$ and $y$ in $M_K,$
 $$ f(1+x)f(1+y)=f(1+xy+x+y).$$
By Equation (\ref{x-exp}) this can be rewritten as
\begin{eqnarray}\label{endo}
\sum_{n=0}^{\infty}a_nx^n \cdot \sum_{n=0}^{\infty}a_ny^n =
\sum_{n=0}^{\infty}a_n(xy+x+y)^n .
\end{eqnarray}
Some computations with Equation (\ref{endo}) give an crucial
identity for $x$:\newline
 For each $m \geq 0,$
\begin{eqnarray}\label{H-der}
a_m\sum_{n=0}^{\infty}a_nx^n =\left(\sum_{n=m}^{\infty}a_n
\binom{n}{m}x^{n-m}\right)(x+1)^m.
\end{eqnarray}
A moment's thought gives that the series in the  parenthesis of
Equation (\ref{H-der}) is nothing but the Hasse derivative of
order $m$ of $f(1+x),$ for which we denote ${\sD}^{(m)}f(1+x)$
(see \cite{H1,HS}). Moreover, the function $f(1+x)$ also appears
on the same equation, so that Equation (\ref{H-der}) can be
rewritten as
\begin{eqnarray}\label{H-der1}
a_mf(1+x)={\sD}^{(m)}f(1+x)(1+x)^m, (m \geq 0).
\end{eqnarray}
Since both $f$ and $(1+x)^m$ take values in $U$ the preceding
identity implies the first part of the result. The second follows
from the first part, as $a_n \in \Fp$ for all $n.$
\end{proof}

Let $f(1+x)=\sum_{n=0}^N a_nx^n \in 1+x\Fp[[x]] $ be a polynomial
of degree $N$ in $\L.$ Since ${\sD}^{(N)}f(1+x)=a_N,$ Equation
(\ref{H-der1}) with $m=N$ gives the following result.

\begin{lemma}\label{ploy}
If a locally  analytic endomorphism
$f(1+x)=\sum_{n=0}^{\infty}a_nx^n \in 1+x\Fp[[x]]$ in $\L$ has
coefficients $a_n =0$ for all but finitely many $n,$ then $f(1+x)
= (1+x)^N$ for some integer $N\geq0.$
\end{lemma}

We next have the same result as Lemma \ref{ploy} for locally
analytic functions $f \in \L$ whose expansion coefficients are
ultimately periodic.

\begin{lemma}\label{rat}
If a locally  analytic endomorphism
$f(1+x)=\sum_{n=0}^{\infty}a_nx^n \in 1+x\Fp[[x]]$ in $\L$ has
coefficients $\{a_n \}_{n\geq 0},$ which are ultimately periodic,
then $f(1+x) = (1+x)^N$ for some integer $N.$
\end{lemma}

\begin{proof}
It is easy to see that if a locally  analytic endomorphism
$f(1+x)=\sum_{n=0}^{\infty}a_nx^n$ has ultimately periodic
coefficients $\{a_n \}_{n\geq 0}$ then $f(1+x)$ is a rational
function in $x.$ Now write $f(1+x)=\frac{P(x)}{Q(x)}$ for two
polynomials $P(x)$ and $Q(x),$ which are relatively prime in
$\Fp[x].$ Since $ f(1)=\frac{P(0)}{Q(0)}=1$ we may take
$P(0)=Q(0)=1,$ if necessary, by cancelling out the same constant.
Also we may here assume that $Q(x)$ is not identically 1 otherwise
we reduce to Lemma \ref{ploy}. Since $f$ is an endomorphism on $U$
then the following equation has to be satisfied:
$$\frac{P(x)}{Q(x)}
\frac{P(y)}{Q(y)} =\frac{P(xy+x+y)}{Q(xy+x+y)}.$$ From this
equation we have
$$P(x)P(y)Q(xy+x+y) =Q(x)Q(y)P(xy+x+y)$$
as polynomials in the polynomial ring $\Fp[x,y]$ in two variables
over $\Fp.$ Since $\Fp[x,y]$ is a unique factorization domain, we
see $P(x)P(y)$ divides $P(xy+x+y)$ and $Q(x)Q(y)$ divides
$Q(xy+x+y).$ The result now follows from Lemma \ref{ploy}, along
with both  $P(x)P(y)=P(xy+x+y)$ and $Q(x)Q(y)=Q(xy+x+y),$ which
come by comparing the highest degrees of polynomials on both sides
of the two preceding equations.
\end{proof}

Before proceeding to prove Theorem 1.1 we mention that it follows
from the general theory of formal groups (see
\cite[Cor.20.2.14]{Haz}) as the formal multiplicative group $G_m$
has height 1, along with Lemma \ref{global} below.  But our goal
in this paper is to give another proof of this result in a
neighborhood of 1 in $U,$ which is very well suited for $G_m$ and
which uses calculus in a profound way. \vskip5pt

\noindent{\bf Proof of Theorem 1.1 } Since $U$ is a module over
$\Zp$ there is a well defined homomorphism from $(\Zp,+)
\rightarrow (\L,\cdot)$ given by $u \mapsto u^y.$ It is relatively
straightforward to check that the map is injective.  Suppose that
$u^y =1$ for some nonzero $y \in \Zp.$ We may here assume that $p$
does not divide $y$ or $y$ is a $p$-adic unit otherwise a $p$th
power mapping reduces to the present case. Since $u$ is
arbitrarily taken  in $U$, consider the expansion of $(1+\pi)^y
=1,$ which is impossible to happen since $y$ is a $p$-adic unit.
It now remains to show the map is surjective, which is the main
point of this note. Take any locally analytic function $f \in \L$
and we may assume by Lemmas \ref{ploy} and \ref{rat} that $f(1+x)
\in 1+x\Fp[[x]]$ is not rational in $x.$ Let $\n_0$ be the maximal
nonnegative integer such that  $f$ is a $p^{\n_{0}}$th power in
$\Fp[[x]].$ Write
$$f(1+x)=f_{0}(1+x)^{p^{\n_{0}}}
\quad {\rm with}\quad f_0 =\sum_{n=0}^{\infty}a_n^{(0)}x^n \in
1+x\Fp[[x]] \quad{\rm in}\quad \L.$$ We here note that $a_1^{(0)}
\not = 0 $ in $\Fp$ otherwise Lemma \ref{ancoeff} forces us to
contradict the maximality of $\n_0.$ Let $g_0(1+x)$ be a locally
analytic endomorphism on $U$ so that $$f_0(1+x)
=(1+x)^{a_1^{(0)}}g_0(1+x).$$ Then we easily check that $g_0(1+x)$
is a locally  analytic endomorphism whose linear coefficient is 0.
Hence by Lemma \ref{ancoeff} there is a maximal positive integer
$\n_1$ such that $g_0(1+x)=f_{1}(1+x)^{p^{\n_{1}}} $ with $f_1
=\sum_{n=0}^{\infty}a_1^{(1)}x^n , a_1^{(1)} \not =0. $ Take $g_1
\in \L $ so that $f_1(1+x) =(1+x)^{a_1^{(1)}}g_1(1+x).$ Then again
by  Lemma \ref{ancoeff} $g_1$ is a $p$th power in $\L.$ Now by
applying the same argument above to $g_1$ and by iterating this
process repeatedly we have two sequences of locally analytic
functions in $\L,$ $\{f_n= \sum_{i=0}^{\infty}a_i^{(n)}x^i,
a_1^{(n)} \not =0 \}_{n\geq0},$ $\{g_n\}_{n\geq0}$ and a sequence
of positive integers $\{\n_n\}_{n\geq1}$ such that
$$f_n(1+x)=(1+x)^{a_1^{(n)}}g_{n}(1+x) \quad {\rm and}\quad
g_{n-1}(1+x)=f_{n}(1+x)^{p^{\n_{n}}} .$$ If we plug in
$\{f_i\}_{0\leq i\leq n}$ and $\{g_i\}_{0 \leq i \leq n}$ into the
function $f$ we have
$$ f(1+x)=(1+x)^{y_n}g_{n}^{p^{\l_n}}$$ with
$\l_n =\n_0+\cdots +\n_n $ and $y_n = a_1^{(0)}p^{\l_0} +
a_1^{(1)}p^{\l_1}+ \cdots + a_1^{(n)}p^{\l_n}.$ Since $\{\n_n
\}_{n\geq1}$ is a sequence of positive integers, as $n \rightarrow
\infty,$ $\l_n \rightarrow \infty$  so $g_{n}^{p^{\l_n}}$
converges to 1. At the same time, $y_n$ converges to some $p$-adic
integer $y.$ Therefore $f(u)=u^y $ for some $y \in \Zp \setminus
\Z.$ The discussion above concludes that if $f$ is a locally
analytic endomorphism on $U,$ then in some open neighborhood of
the identity $f$ is  of the form $f(1+x)=(1+x)^y$ for some
$p$-adic integer. By Lemma \ref{global} below it follows that the
map is surjective. $\Box$

\begin{lemma}\label{global}
If a locally analytic endomorphism $f$ on $U$ is of the form
$f(1+x)= (1+x)^y$ for some $y \in \Zp,$ in some open neighborhood
of 1, then it is identically  $(1+x)^y.$
\end{lemma}
\begin{proof}
Let $u=1+x$ be an arbitrary element of $U$ and let $V$ be some
open neighborhood of 1.
 Then we see for some positive integer $n$ that
$(1+x)^{p^n}= 1+ x^{p^n}$ is in $V.$ Thus we must have
$${f(1+x)}^{p^n} =f((1+x)^{p^n})=f(1+{x}^{p^n})=(1+x^{p^n})^y
=((1+x)^{y})^{p^n},$$ for some $p$-adic integer $y.$ Therefore the
injection of $p$th power mappings gives $f(1+x)=(1+x)^y$ on $U.$

\end{proof}
\vskip5pt

\noindent{\bf Proof of Corollary 1.2 } We note that
$f(1+x)=\sum_{n=0}^{\infty}a_nx^n \in 1+x\Fp[[x]]$ in $\L$ is a
locally analytic automorphism on $U$ if and only if $a_1 \not =0.$
We here leave a justification of this assertion to the reader.
From the proof of Theorem \ref{main} we have $f(u)=u^y $ for some
$p$-adic unit $y.$  So we are done. $\Box$

\vskip5pt
Now the question arises of whether or not exponential
functions $u^y (y \in \Z_p\setminus\Z)$ have ultimately periodic
coefficients. To this end, we first state the well known result
for the $p$-adic numbers.

\begin{proposition}\label{pres}
A $p$-adic number $a=\sum_{n=\w \in \Z}^{\infty}a_np^n$ $ (0\leq
a_i<p)$ is rational in $\Q$ if and only if the sequence of digits
$\{a_n \}$ is ultimately periodic.
\end{proposition}
\begin{proof}
See \cite[p. 147]{H2}.
\end{proof}

We here have a function field analogue of Proposition \ref{pres}.

\begin{proposition}\label{ffperiodic}
Let $\Fq((x))$ be the field of formal Laurent series in one
variable $x$ over a finite field $\Fq.$ Then $f(x)=\sum_{n=\w \in
\Z}^{\infty}a_nx^n$ is a rational function in $x$ if and only if
the sequence of elements $\{a_n \} \subset \Fq$ is ultimately
periodic.
\end{proposition}
\begin{proof}
The translation to $\Fq((x))$ of the arguments\cite[p. 147]{H2} in
the proof of Proposition \ref{pres} goes in a parallel way as in
the classical case. For it must rely on the fact that $\Fq$ is
finite, hence we have the well known analogue of Euler's theorem
(see \cite[Prop.1.8 p. 5]{R}).
\end{proof}

\vskip5pt
\noindent{\bf Proof of Corollary 1.3 } The equivalence
of (1) and (2) follows from Proposition \ref{ffperiodic}. And we
see that the equivalence of (1) and (3) follows from the proof of
Theorem \ref{main}.$\Box$

\vskip5pt
We close the paper with several remarks.

\noindent 1.  The methodology of Theorem 1.1 reminds us of the
well known result in calculus that every differentiable function
$y$ with $y^\prime =y$ must be a multiple of $e^x.$ We take one
such $h$ and divide by $e^x$ and then show the derivative of the
resulting function is identically 0. We can also use the Hasse
derivatives or divided derivatives to calculate something similar
in characteristic $p.$

\noindent 2.  It is straightforward to see that the result for
$\Zp$ analogous to Theorem 1.1 is trivial. But it might be
nontrivial to have some parallel results on the ring of integers
of finite extensions of $\Qp.$

\noindent 3. $(1) \Leftrightarrow (3)$ in Corollary 1.3 also
follows from the known result by Mendes France and van der Poorten
\cite{MvP} using the finite automata tool: $(1+x)^y$ is algebraic
over $\Fq(x)$ if and only if $y$  is in $\Zp \cap \Q.$ There their
result also holds if $1+x$ is replaced with an algebraic formal
power series with constant term equal to 1.


{\centerline {\bf Acknowledgements}

The author thanks David Goss for his invaluable comments on the
earlier draft of the paper and for forwarding the paper \cite{MvP}
to the author via the email sent by J. P. Allouche.

\bibliographystyle{amsplain}

\noindent{Department of Mathematics, Inha University, Incheon,
Korea 402-751} \newline \noindent{E-mail address: \em
stj@inha.ac.kr}

\end{document}